\documentclass[a4paper,11pt]{amsart}

\usepackage{amsmath,amsthm,amssymb,enumerate}
\usepackage{epsfig}
\usepackage{amssymb}
\usepackage{amsmath}
\usepackage{amssymb}
\usepackage{amsmath,amsthm}
\usepackage[latin1]{inputenc}
\usepackage[T1]{fontenc}
\usepackage{path}
\usepackage{ae,aecompl}
\usepackage{amsfonts}
\usepackage{amsxtra}
\usepackage{euscript,mathrsfs}
\usepackage{color}
\usepackage[left=2.5cm,right=2.5cm,top=4cm,bottom=4cm]{geometry}
\usepackage[citecolor=blue,colorlinks=true]{hyperref}
\allowdisplaybreaks

\usepackage{enumitem}
\setenumerate{label={\rm (\alph{*})}}

\usepackage{amsfonts}
\usepackage{amsxtra}

\usepackage[frak,hyperref,theorem_section]{paper_diening}
\numberwithin{equation}{section}

\newcommand{\Div}{\divergence}
\newcommand{\ep}{\bfvarepsilon}
\newcommand{\R}{\mathbb R}

\newcommand{\E}{\mathbb E}
\newcommand{\p}{\mathbb P}

\newcommand{\F}{\mathscr F}

\newcommand{\dd}{\mathrm d}

\newcommand{\dt}{\, \mathrm{d}t}

\newcommand{\dif}{\mathrm{d}}
\newcommand{\totdif}{\mathrm{D}}
\newcommand{\mf}{\mathscr{F}}

\begin{document}

\title[Time regularity for SPDE\MakeLowercase{s}]{On time regularity of stochastic evolution equations with monotone coefficients}

\author{Dominic Breit}
\address[D. Breit]{Department of Mathematics, Heriot-Watt University, Riccarton Edinburgh EH14 4AS, UK}
\email{d.breit@hw.ac.uk}

\author{Martina Hofmanov\'a}
\address[M. Hofmanov\'a]{Technical University Berlin, Institute of Mathematics, Stra\ss e des 17. Juni 136, 10623 Berlin, Germany}
\email{hofmanov@math.tu-berlin.de}

\begin{abstract}
We report on a time regularity result for stochastic evolutionary PDEs with monotone coefficients. If the diffusion coefficient is bounded in time without additional space regularity we obtain a fractional Sobolev type time regularity of order up to $\tfrac{1}{2}$ for a certain functional $ G( u )$ of the solution. Namely, $ G( u )=\nabla u $ in the case of the heat equation and $G( u )=|\nabla u |^{\frac{p-2}{2}}\nabla u $ for the $p$-Laplacian. The motivation is twofold. On the one hand, it turns out that this is the natural time regularity result that allows to establish the optimal rates of convergence for numerical schemes based on a time discretization. On the other hand, in the linear case, i.e. where the solution is given by a stochastic convolution, our result complements the known stochastic maximal space-time regularity results for the borderline case not covered by other methods.
\end{abstract}

\subjclass[2010]{}

\date{\today}

\maketitle

%
%
%
%
%
%
%
%
%
%


\section{Time regularity}

Let $H,\,U$ be separable Hilbert spaces and let $V$ be a Banach space such that $ V \hookrightarrow  H \hookrightarrow  V '$ is a Gelfand triple with continuous and dense embeddings. We are interested in stochastic evolution equations of the form
\begin{align}\label{eq:}
\begin{aligned}
\dd u&= A  (t, u )\dt+B(t,u)\,\dd W,\\
u(0)&=u_0,
\end{aligned}
\end{align}
where $W$ is a $U$-valued cylindrical Wiener process on a probability space $(\Omega,\F,\p)$ with a normal filtration $(\mf_t)$ and the maps
$$A:\Omega\times[0,T]\times V\to V',\qquad B:\Omega\times[0,T]\times H\to L_2(U;H)$$
are $(\mf_t)$-progressively measurable and satisfy
\begin{enumerate}[label={\rm (H\arabic{*})}, leftmargin=*]
\item[(H1)]\label{eq:mu2'} Monotonicity: there exists $c_1\in\R$ such that for all $ u , v \in  V $, $t\in[0,T]$
\begin{align*}
2{}_{ V '}\langle  A  (t, u )- A  (t,v ), u - v \rangle_{ V }+\|B(t,u)-B(t,v)\|_{L_2(U;H)}^2\leq c_1\|u-v\|_H^2.
\end{align*}
\item[(H2)]\label{eq:mu0} Hemicontinuity: for all $ u , v , w   \in  V $, $\omega\in\Omega$ and $t\in[0,T],$ the mapping
$$\R\ni \lambda\mapsto{}_{V'}\langle A  (\omega,t, u +\lambda v ), w   \rangle_V$$
is continuous.
\item[(H3)]\label{eq:mu} Coercivity: there exist $q\in(1,\infty)$, $c_2\in[0,\infty)$, $c_3\in\R $ such that for all $ u \in V$, $t\in[0,T]$
\begin{align*}
{}_{ V '}\langle  A  (t, u ), u \rangle_{ V }\leq -c_2\| u \|_V^q+c_3.
\end{align*}
\item[(H4)]\label{eq:mu2} Growth of $A$: there exists $c_4\in (0,\infty)$ such that for all $ u \in  V $, $t\in[0,T]$
\begin{align*}
\| A  (t, u )\|_{ V '}^{q'}\leq c_4 \big(1+\| u \|^q_{ V }\big).
\end{align*}
\item[(H5)] Growth of $B$: there exists $c_5\in (0,\infty)$ and $(\mathscr F_t)$-adapted $f\in L^2(\Omega;L^\infty(0,T))$ such that for all $ u \in  H $, $t\in[0,T]$
\begin{align*}
\|B(t,u)\|_{L_2(U;H)}\leq c_5 (f+\|u\|_H).
\end{align*}
\end{enumerate}

The literature devoted to the study of these equations is quite extensive. The question of existence of a unique (variational) solution to equations of the form \eqref{eq:} is well-understood: first results were established in \cite{Pa,KrRo2}, for an overview in the above stated generality and further references we refer the reader to \cite{PrRo}. Existence of a strong solution under various assumptions appeared in \cite{Br3,G} and numerical approximations were studied in \cite{GM1,GM2}.
In the case of linear operator $A$ which generates a strongly continuous semigroup, more is known concerning regularity and maximal regularity (see e.g. \cite{B,Ho2, vNVW}).

Naturally, the time regularity of a solution to \eqref{eq:} is limited by the regularity of the driving Wiener process $W$. In particular, since the trajectories of $W$ are only $\alpha$-H\"older continuous for $\alpha<\tfrac{1}{2}$, it can be seen from the integral formulation of \eqref{eq:} that the trajectories of $ u $ are $\alpha$-H\"older continuous as functions taking values in $ V '$.
This can be improved if some additional regularity in space of the solution is known, that is, the equation is satisfied in a stronger sense.
In this note, we are particularly interested in situations where such additional space regularity is either not available or limited. This is typically the case when
\begin{itemize}
\item[(i)] $ A  $ is linear but the noise is not smooth enough: if $u$ is a variational solution to \eqref{eq:} then the standard assumption is $B(u)\in L^2_{w^*}(\Omega;L^\infty(0,T;L_2( U ; H ))$.\footnote{Here $L^2_{w^*}(\Omega;L^\infty(0,T;L_2(U;H))$ is the space of weak$^*$-measurable mappings $h:\Omega\to L^\infty(0,T;L_2(U;H))$ such that $\E\esssup_{0\leq t\leq T}\|h\|^2_{L_2(U;H)}<\infty.$}
 \item[(ii)] $ A  $ is nonlinear as for instance the $p$-Laplacian $ A  ( u )=\Div(|\nabla u |^{p-2}\nabla u )$ or a more general nonlinear operator with $p$-growth and, in addition, the noise represents the same difficulty as in (i).
\end{itemize}

In order to formulate our main result we need several additional assumptions upon the operator $A$ and the initial datum $u_0$. On the one hand, we introduce a notion of $G$-monotonicity which represents a stronger version of the monotonicity assumption on $A$, on the other hand, we suppose certain regularity in time of $A$ as well as regularity of the initial condition. To be more precise, we assume
\begin{itemize}
\item[(H6)]\label{eq:gmon} $G$-monotonicity: there exists a bounded (possibly nonlinear) mapping $G: V\to H$ and $c_6\in (0,\infty)$ such that for all $u,v\in V$, $t\in[0,T]$
\begin{align*}
-{}_{ V '}\langle  A  (t, u )- A  (t, v ), u - v \rangle_{ V }\geq c_6 \| G( u )- G( v )\|_{ H }^2.
\end{align*}
\item[(H7)] Time regularity of $A$: there exists $c_7\in (0,\infty)$ such that for all $u\in V$, $t,s\in[0,T]$
\begin{align*}
\|A(t,u)-A(s,u)\|^{q'}_{V'}\leq c_7\big(\|u\|^q_V+1\big)|t-s|.
\end{align*}
\item[(H8)] Regularity of $u_0$: $ A (t, u_0)\in H$ a.s. for all $t\in[0,T]$ and there exists $c_8\in (0,\infty)$ such that
\begin{align*}
\sup_{0\leq t\leq T}\E\|A(t,u_0)\|_H^2\leq c_8.
\end{align*}
\end{itemize}


Note that it can be readily checked that the operators $A$ in the above mentioned examples (i) and (ii) are $G$-monotone. Indeed, if $ A  $ is linear and symmetric negative definite we can choose $ G = (-A)  ^{1/2}$ and, as was shown in \cite{DieE08}, the $p$-Laplacian is covered via $ G ( u )=|\nabla u |^{\frac{p-2}{2}}\nabla u $ which is the natural quantity to establish its regularity properties.


Finally we have all in hand to state our result.
\begin{theorem}\label{eq:time}
Assume that {\em(H1)-(H8)} hold true. If $ u $ is a solution to \eqref{eq:}, in particular
\begin{align}\label{eq:210}
u\in L^q(\Omega;L^q(0,T;V))\cap L^2_{w^*}(\Omega;L^\infty(0,T;H)),
\end{align}
then
\begin{align}\label{eq:fract}
 G ( u )\in L^2(\Omega;W^{\alpha,2}(0,T;H))\quad\text{ for all }\quad\alpha<\tfrac{1}{2}.
\end{align}
\end{theorem}

\begin{remark}
If one drops the assumption (H8) then \eqref{eq:fract} holds locally in time away from 0. 
\end{remark}

\begin{corollary}\label{rem:new}
The statement of Theorem \ref{eq:time} continues to hold if we replace {\em (H6)} with the following assumption:
\begin{itemize}
\item[{\em (H6')}]\label{eq:gmon'} modified $G$-monotonicity: there exists a separable Hilbert space $\mathcal H$ (generally different from $H$) and a bounded mapping $G: V\to \mathcal H$ and $c_6'\in (0,\infty)$ such that for all $u,v\in V$, $t\in[0,T]$
\begin{align*}
-{}_{ V '}\langle  A  (t, u )- A  (t, v ), u - v \rangle_{ V }\geq c_6' \| G( u )- G( v )\|_{\mathcal H }^2.
\end{align*}
\end{itemize}
In this case we have to replace $H$ by $\mathcal H$ in \eqref{eq:210} and \eqref{eq:fract}.
\end{corollary}

Let us now explain what are the main motivations for such a result.
First, it turns out that \eqref{eq:fract} is the natural time regularity that allows to establish the optimal rates of convergence for numerical schemes based on time discretization (or a space-time discretization provided a suitable space regularity can be proved as well). Indeed, with this time regularity at hand, a finite element based space-time discretization of stochastic $p$-Laplace type systems will be studied in \cite{BHLL}.
 A similar strategy can be directly applied to establish rates of convergence for time discretization of more general monotone SPDEs satisfying (among others) the key $G$-monotonicity assumption.

Second, if $A$ is a linear infinitesimal generator of a strongly continuous semigroup $S$ on $H$ then the (mild) solution to \eqref{eq:} with $u_0=0$ is given by the stochastic convolution
$$u(t)=\int_0^t S(t-s)B(s,u_s)\,\dif W_s$$
and our result gives $u\in L^2(\Omega;W^{\alpha,2}(0,T;D((-A)^{1/2}))$. Recall that the space $D((-A)^{1/2})$ here is the borderline case regarding regularity for the stochastic convolution, namely, $(-A)^{1/2}u$ may not even have a pathwise continuous version whereas for $(-A)^{1/2-\varepsilon}u$ has $\alpha$-H\"older continuous trajectories for $\alpha\in (0,\varepsilon)$ (see \cite[Theorem 5.16, Subsection 5.4.2]{PrZa}). Consequently, the borderline case is typically not covered by known methods such as factorization \cite{PrZa, B} or stochastic maximal regularity (see \cite[Theorem 1.1, Theorem 1.2]{vNVW}) and Theorem \ref{eq:time} provides an additional information based on a rather simple argument.

\begin{proof}[Main ideas of the proof of Theorem \ref{eq:time}:]
A complete proof will be given in \cite{BHLL}. It is based on a new version of the It\^o formula which applies to time differences and yields the following:
let $0<h\ll1$ and $t\in(h,T]$ then it holds true a.s.
\begin{equation}\label{eq:ito}
\begin{split}
\|u(t)-u({t-h})\|_H^2&=\|u(h)-u_0\|_H^2+2\int_h^t{}_{V}\langle u({\sigma})-u({\sigma-h}),\dd u({\sigma})\rangle_{V'}\\
&\quad-2\int_0^{t-h}{}_{V}\langle u({\sigma+h})-u({\sigma}),\hat{\dd} u({\sigma})\rangle_{V'}+\langle\!\langle u\rangle\!\rangle_t-\langle\!\langle u\rangle\!\rangle_h-\langle\!\langle u\rangle\!\rangle_{t-h}.
\end{split}
\end{equation}
Here $\hat{\dd}u$ denotes the backward It\^o stochastic differential and $\langle\!\langle \cdot\rangle\!\rangle$ the quadratic variation process. The appearance of the backward It\^o stochastic integral comes from the fact that the It\^o formula is applied to the time difference $t\mapsto u(t)-u(t-h)$. Indeed, if $M$ denotes the martingale part of $u$, then
for every fixed $t_0\in[0,T)$ the process $t\mapsto M_t-M_{t_0}$
is a (forward) local martingale with respect to the forward filtration given by
$\sigma(M_r-M_{t_0};\,{t_0}\leq r\leq t)$, $t\in[{t_0},T),$
whereas for every fixed $t_1\in[0,T]$ the process $t\mapsto M_{t_1}-M_t$ is a (backward) local martingale with respect to the backward filtration given by
$\sigma(M_{t_1}-M_r;\,t\leq r\leq t_1)$, $t\in[0,t_1].$

As the next step, we substitute for $\dd u$ and $\hat{\dd} u$ in \eqref{eq:ito}, take expectation and apply hypotheses (H5)-(H8). Finally we obtain that
\begin{align*}
\frac{1}{h}\,\E\int_0^{T-h} \big\|G(u(\sigma+h))-G(u(\sigma))\big\|_{H}^2\,\dd \sigma&\leq C
\end{align*}
which implies the required regularity.
%
%
%
\end{proof}

\section{Applications}

In this section we present some concrete examples of problems which are covered by our result.

\subsection{The linear case}
Let us assume that $A:D(A)\subset H\to H$ is linear dissipative and symmetric infinitesimal generator of a strongly continuous semigroup on $H$. Then the square root $(-A)^{1/2}$ is well-defined and setting $V=D((-A)^{1/2})$ (equipped with the graph norm) we obtain, for all $u,v\in V$, that
\begin{align*}
-\langle A u-A v,u- v\rangle_{H}=\big\|(-A)^{1/2}u-(-A)^{1/2}v\big\|_{ H}^2=\|u-v\|_{V}^2.
\end{align*}
Thus the hypothesis (H6) holds true with $ G=(-A)^{1/2}$ and Theorem \ref{eq:time} applies.

\subsection{The $p$-Laplace type systems}\label{subsec:plapla}
Let $\mathcal O\subset\R^d$ be a bounded Lipschitz domain and let $H=L^2(\mathcal O)$. We suppose that $\Phi$ satisfies (H1) and (H5). We are interested in the system
\begin{align*}
\begin{aligned}
\dd \bfu &=\Div \bfS(\nabla \bfu )\dt+\Phi( \bfu )\dd W,\\
\bfu|_{\partial\mathcal O}&=0,\\
\bfu(0)&=\bfu_0,
\end{aligned}
\end{align*}
where $\bfS:\R^{d\times D}\rightarrow \R^{d\times D}$ is a general nonlinear operator with $p$-growth, i.e.
\begin{align*}
c (\kappa+|\bfxi|)^{p-2}|\bfzeta|^2\leq \totdif\bfS(\bfxi)(\bfzeta,\bfzeta)\leq C(\kappa+|\bfxi|)^{p-2} |\bfzeta|^2
\end{align*}
for all $\bfxi,\bfzeta\in\R^{d\times D}$ with some constants $c,\,C>0$, $\kappa\geq0$ and $p\in[\frac{2d}{d+2},\infty)$. Then the assumptions (H1)-(H4) are satisfied with $V=W^{1,p}_0(\mathcal O)$ and, in addition, it is well known from the deterministic setting (and was already discussed in \cite{Br3} in the stochastic setting) that an important role for this system is played by the function
$$\bfF(\bfxi)=(\kappa+|\bfxi|)^{\frac{p-2}{2}}\bfxi.$$
It is used in regularity theory \cite{AF} and also for the numerical approximation \cite{BaLi1,DiRu}.
The essential property of $\bfF$ can be characterized by the inequality
\begin{align*}
\lambda|\bfF(\bfxi)-\bfF(\bfeta)|^2\leq \big(\bfS(\bfxi)-\bfS(\bfeta)\big):(\bfxi-\bfeta)\leq \Lambda|\bfF(\bfxi)-\bfF(\bfeta)|^2\qquad\forall\bfxi,\bfeta\in\R^{d\times D}
\end{align*}
for some positive constants $\lambda,\,\Lambda$ depending only on $p$ (see for instance \cite{DieE08}).
Consequently, for all $\bfu , \,\bfv\in V$,
\begin{align*}
\lambda \big\|\bfF(\nabla \bfu )-\bfF(\nabla \bfv)\big\|_{H}^2\leq -{}_{V'}\langle \Div \bfS(\nabla \bfu )-\Div \bfS(\nabla \bfv), \bfu -\bfv\rangle_{V}\leq \Lambda \big\| \bfF(\nabla \bfu )- \bfF(\nabla \bfv)\big\|_{H}^2,
\end{align*}
and therefore (H6) is satisfied and Theorem \ref{eq:time} yields
\begin{align*}
\bfF(\nabla \bfu )\in L^2(\Omega;W^{\alpha,2}(0,T;L^2(\mathcal O)))\quad\text{ for all }\quad\alpha<\tfrac{1}{2}.
\end{align*}
Note that in case of the heat equation (i.e. $p=2$)
the operator $F$ is the identity.

\subsection{The $p$-Stokes system}

In continuum mechanics, the motion of a homogeneous incompressible fluid is described by its velocity field $\bfu $ and its pressure function $\pi$. If the flow is slow motion can be described via the system
\begin{align}\label{eq:stokes}
\begin{aligned}
\dd \bfu& =\Div \bfS(\ep( \bfu ))\dt+\nabla\pi\dt+\Phi( \bfu )\dd W,\\
\Div  \bfu &=0,\\
 \bfu|_{\partial \mathcal O} &=0,\\
 \bfu(0)&=\bfu_0,
 \end{aligned}
\end{align}
where $\mathcal O$ and $\bfS$ satisfy the hypotheses of Subsection \ref{subsec:plapla} and $\ep(\bfu)=\tfrac{1}{2}\big(\nabla\bfu+\nabla\bfu^T\big)$ is the symmetric gradient of the velocity field $\bfu$.
In comparison to the Navier--Stokes system the convective term $-(\nabla \bfv) \bfv\dt$ on the right-hand-side of the momentum equation \eqref{eq:stokes}$_1$ is neglected (see \cite{Br2} for the corresponding Navier--Stokes system for power-law fluids and further references). In the following functional analytical setting 
\begin{align*}
H&=L^2_{\Div}(\mathcal O)=\overline{C^\infty_{0,\Div}(\mathcal O)}^{L^2(\mathcal O)},\quad
 V=W^{1,p}_{0,\Div}(\mathcal O)=\overline{C^\infty_{0,\Div}(\mathcal O)}^{W^{1,p}(\mathcal O)}
\end{align*}
where
$$C^\infty_{0,\Div}(\mathcal O)=\{ \bfw\in C^\infty_0(\Omega):\,\,\Div \bfw=0\},$$
the pressure function does not appear.
Similarly to the $p$-Laplace system we set
 $G( \bfu )= \bfF(\ep( \bfu ))$ and obtain, for all $\bfu ,\,\bfv\in V$,
\begin{align*}
 \lambda \big\| \bfF(\ep(\bfu ))- \bfF(\ep( \bfv))\big\|_{\mathcal H}^2\leq - {}_{V'}\langle \Div \bfS(\nabla \bfu )-\Div \bfS(\nabla \bfv), \bfu - \bfv\rangle_{V}\leq \Lambda \big\| \bfF(\ep( \bfu ))- \bfF(\ep( \bfv))\big\|_{\mathcal H}^2,
\end{align*}
where $\mathcal H=L^2(\mathcal O)$. Corollary \ref{rem:new} applies and we gain
\begin{align*}
\bfF(\ep(\bfu))\in L^2(\Omega;W^{\alpha,2}(0,T;L^2(\mathcal O)))\quad\text{ for all }\quad\alpha<\tfrac{1}{2}.
\end{align*}


\end{document}